\documentclass[11pt,leqno]{article}

\usepackage{amsmath,amsthm,amsfonts}
\usepackage{amssymb}
\usepackage{euscript}
\usepackage[all]{xy}

\theoremstyle{plain}

\theoremstyle{definition}

\newcommand{\acat}[1]{\mathsf{a}_{#1}}

\newcommand{\lcat}{\mathsf{l}}
\newcommand{\sets}{\mathsf{Sets}}

\newcommand{\sfam}{\mathcal{M}}

\newcommand{\defs}{\mathsf{Def}_{\sfam}}
\newcommand{\defm}[1]{\mathsf{Def}_{#1}}

\newcommand{\ch}{\mathsf{YH}}
\newcommand{\cpd}{*}

\newcommand{\z}{\mathbb{Z}}

\newcommand{\aff}{\mathbf{A}}

\newcommand{\diff}[1]{\mathrm{Diff}({#1})}

\DeclareMathOperator{\enm}{End}
\DeclareMathOperator{\ext}{Ext}

\DeclareMathOperator{\hmm}{Hom}
\DeclareMathOperator{\id}{id}

\DeclareMathOperator*{\osum}{\oplus}

\DeclareMathOperator{\spec}{Spec}

\begin{document}

\title{\LARGE\bf An example of noncommutative deformations}
\author{\Large Eivind Eriksen \\ \\ Oslo University College \\
Postboks 4, St. Olavs plass \\ N-0130 Oslo, Norway \\ \\ E-mail:
eeriksen@hio.no}
\date{\today}

\maketitle

\thispagestyle{empty}

\begin{abstract}
We compute the noncommutative deformations of a family of modules over
the first Weyl algebra. This example shows some important properties of
noncommutative deformation theory that separates it from commutative
deformation theory.
\end{abstract}

\section{Introduction}

Let $k$ be an algebraically closed field and let $A$ be an associative
$k$-algebra. For any left $A$-module $M$, there is a flat commutative
deformation functor
    \[ \defm M: \lcat \to \sets \]
defined on the category $\lcat$ of local Artinan commutative
$k$-algebras with residue field $k$. We recall that for an object $R
\in \lcat$, a flat deformation of $M$ over $R$ is a pair $(M_R, \tau)$,
where $M_R$ is an $A$-$R$ bimodule (on which $k$ acts centrally) that
is $R$-flat, and $\tau: M_R \otimes_R k \to M$ is an isomorphism of
left $A$-modules. Moreover, $(M_R, \tau) \sim (M'_R, \tau')$ as
deformations in $\defm M(R)$ if there is an isomorphism $\eta: M_R \to
M'_R$ of $A$-$R$ bimodules such that $\tau = \tau' \circ (\eta \otimes
1)$.

Laudal introduced noncommutative deformations of modules in Laudal
\cite{laud02}. For any finite family $\sfam = \{ M_1, \dots, M_p \}$ of
left $A$-modules, there is a noncommutative deformation functor
    \[ \defs: \acat p \to \sets \]
defined on the category $\acat p$ of $p$-pointed Artinian $k$-algebras.
We recall that an object $R$ of $\acat p$ is an Artinian ring $R$,
together with a pair of structural ring homomorphisms $f: k^p \to R$
and $g: R \to k^p$, such that $g \circ f = \id$ and the radical $I(R) =
\ker(g)$ is nilpotent. The morphisms of $\acat p$ are ring
homomorphisms that commute with the structural morphisms.

A deformation of the family $\sfam$ over $R$ is a $(p+1)$-tuple $(M_R,
\tau_1, \dots, \tau_p)$, where $M_R$ is an $A$-$R$ bimodule (on which
$k$ acts centrally) such that $M_R \cong ( M_i \otimes_k R_{ij} )$ as
right $R$-modules, and $\tau_i: M_R \otimes_R k_i \to M_i$ is an
isomorphism of left $A$-modules for $1 \le i \le p$. By definition,
    \[ ( M_i \otimes_k R_{ij} ) = \osum_{1 \le i,j \le p} \;
    M_i \otimes_k R_{ij} \]
with the natural right $R$-module structure, and $k_1, \dots, k_p$ are
the simple left $R$-modules of dimension one over $k$. Moreover, $(M_R,
\tau_1, \dots, \tau_p) \sim (M'_R, \tau'_1, \dots, \tau'_p)$ as
deformations in $\defs(R)$ if there is an isomorphism $\eta: M_R \to
M'_R$ of $A$-$R$ bimodules such that $\tau_i = \tau'_i \circ (\eta
\otimes 1)$ for $1 \le i \le p$.

There is a cohomology theory and an obstruction calculus for $\defs$,
see Laudal \cite{laud02} and Eriksen \cite{erik06}. We compute the
noncommutative deformations of a family $\sfam = \{ M_1, M_2 \}$ of
modules over the first Weyl algebra using the constructive methods
described in Eriksen \cite{erik06}.

\section{An example of noncommutative deformations of a family}

Let $k$ be an algebraically closed field of characteristic $0$, let $A
= k[t]$, and let $D = \diff A$ be the first Weyl algebra over $k$. We
recall that $D = k\langle t, \partial\rangle / ( \partial \; t - t \;
\partial - 1)$. Let us consider the family $\sfam = \{ M_1, M_2 \}$ of
left $D$-modules, where $M_1 = D/D \cdot \partial \cong A$ and $M_2 =
D/D \cdot t \cong k[\partial]$. We shall compute the noncommutative
deformations of the family $\sfam$.

In this example, we use the methods described in Eriksen \cite{erik06}
to compute noncommutative deformations. In particular, we use the
cohomology $\ch^n(M_j,M_i)$ of the Yoneda complex
    \[ YC^p(M_j,M_i) = \prod_{m \ge 0} \quad \hmm_D(L_{m,j},L_{m-p,i}) \]
for $1 \le i,j \le 2$, where $(L_{\cpd,i},d_{\cpd,i})$ is a free
resolution of $M_i$, and an obstruction calculus based on these free
resolutions. We recall that $\ch^n(M_j,M_i) \cong \ext^n_D(M_j,M_i)$.

Let us compute the cohomology $\ch^n(M_j,M_i)$ for $n = 1,2$, $1 \le
i,j \le 2$. We use the free resolutions of $M_1$ and $M_2$ as left
$D$-modules given by
\begin{align*}
0 \gets & M_1 \gets D \xleftarrow{\cdot \partial} D \gets 0 \\
0 \gets & M_2 \gets D \xleftarrow{\cdot t} D \gets 0
\end{align*}
and the definition of the differentials $YC^0(M_j,M_i) \to
YC^1(M_j,M_i) \to YC^2(M_j,M_i) = 0$ in the Yoneda complex, and obtain
\begin{align*}
    \ch^1(M_1,M_1) \cong \ext^1_D(M_1,M_1) &= 0 &
    \ch^1(M_1,M_2) \cong \ext^1_D(M_1,M_2) &= k \cdot \xi_{21} \\
    \ch^1(M_2,M_1) \cong \ext^1_D(M_2,M_1) &= k \cdot \xi_{12} &
    \ch^1(M_2,M_2) \cong \ext^1_D(M_2,M_2) &= 0
\end{align*}
The base vector $\xi_{ij}$ is represented by the $1$-cocycle given by
$D \xrightarrow{\cdot 1} D$ in $YC^1(M_j,M_i)$ when $i \neq j$. Since
$YC^2(M_j,M_i) = 0$ for all $i,j$, it is clear that $\ch^2(M_j,M_i)
\cong \ext^2_D(M_j,M_i) = 0$ for $1 \le i,j \le 2$.

We conclude that $\defm \sfam$ is unobstructed. Hence, in the notation
of Eriksen \cite{erik06}, the pro-representing hull $H$ of $\defs$ is
given by
    \[ H = \left( \begin{matrix} H_{11} & H_{12} \\ H_{21} & H_{22}
    \end{matrix} \right) \cong \left( \begin{matrix} k[[ s_{12} s_{21}
    ]] & \langle s_{12} \rangle \\ \langle s_{21} \rangle & k[[ s_{21}
    s_{12} ]] \end{matrix} \right) \]
where $\langle s_{12} \rangle = H_{11} \cdot s_{12} \cdot H_{22}$ and
$\langle s_{21} \rangle = H_{22} \cdot s_{21} \cdot H_{11}$.

In order to describe the versal family $\sfam_H$ of left $D$-modules
defined over $H$, we use M-free resolutions in the notation of Eriksen
\cite{erik06}. In fact, the $D$-$H$ bimodule $\sfam_H$ has an M-free
resolution of the form
    \[ 0 \gets \sfam_H \gets \left( \begin{matrix} D
    \widehat{\otimes}_k H_{11} & D \widehat{\otimes}_k H_{12} \\ D
    \widehat{\otimes}_k H_{21} & D \widehat{\otimes}_k H_{22}
    \end{matrix} \right) \xleftarrow{d^H} \left( \begin{matrix} D
    \widehat{\otimes}_k H_{11} & D \widehat{\otimes}_k H_{12} \\ D
    \widehat{\otimes}_k H_{21} & D \widehat{\otimes}_k H_{22}
    \end{matrix} \right) \gets 0 \]
where $d^H = (\cdot \partial) \widehat{\otimes} e_i - (\cdot 1)
\widehat{\otimes} s_{12} - (\cdot 1) \widehat{\otimes} s_{21} + (\cdot
t) \widehat{\otimes} e_2$. This means that for any $P,Q \in D$, we have
that $d^H(P \otimes e_1) = (P \cdot \partial) \widehat{\otimes} e_1 -
(P \cdot 1) \widehat{\otimes} s_{21}$ and $d^H(Q \otimes e_2) = (Q
\cdot t) \widehat{\otimes} e_2 - (Q \cdot 1) \widehat{\otimes} s_{12}$.

We remark that there is a natural algebraization $S$ of the
pro-representation hull $H$, given by
    \[ S = \left( \begin{matrix} S_{11} & S_{12} \\ S_{21} & S_{22}
    \end{matrix} \right) \cong \left( \begin{matrix} k[ s_{12} s_{21}
    ] & \langle s_{12} \rangle \\ \langle s_{21} \rangle & k[ s_{21}
    s_{12} ] \end{matrix} \right) \]
In other words, $S$ is an associative $k$-algebra of finite type such
that the $J$-adic completion $\widehat S \cong H$ for the ideal $J =
(s_{12}, s_{21}) \subseteq S$. The corresponding algebraization
$\sfam_S$ of the versal family $\sfam_H$ is given by the M-free
resolution
    \[ 0 \gets \sfam_S \gets \left( \begin{matrix} D \otimes_k S_{11}
    & D \otimes_k S_{12} \\ D \otimes_k S_{21} & D \otimes_k S_{22}
    \end{matrix} \right) \xleftarrow{d^S} \left( \begin{matrix} D
    \otimes_k S_{11} & D \otimes_k S_{12} \\ D \otimes_k S_{21} & D
    \otimes_k S_{22} \end{matrix} \right) \gets 0 \]
with differential $d^S = (\cdot \partial) \otimes e_i - (\cdot 1)
\otimes s_{12} - (\cdot 1) \otimes s_{21} + (\cdot t) \otimes e_2$.

We shall determine the $D$-modules parameterized by the family
$\sfam_S$ over the noncommutative algebra $S$ --- this is much more
complicated than in the commutative case. We consider the simple left
$S$-modules as the points of the noncommutative algebra $S$, following
Laudal \cite{laud03}, \cite{laud05}. For any simple $S$-module $T$, we
obtain a left $D$-module
    \[ M_T = \sfam_S \otimes_S T \]
Therefore, we consider the problem of classifying simple $S$-modules of
dimension $n \ge 1$.

Any $S$-module of dimension $n \ge 1$ is given by a ring homomorphism
$\rho: S \to \enm_k(T)$, and we may identify $\enm_k(T) \cong M_n(k)$
by choosing a $k$-linear base $\{ v_1, \dots, v_n \}$ for $T$. We see
that $S$ is generated by $e_1, s_{12}, s_{21}$ as a $k$-algebra (since
$e_2 = 1 - e_1$), and there are relations
    \[ s_{12}^2 = s_{21}^2 = 0, \; e_1^2 = e_1, \; e_1 s_{12} = s_{12},
    \; s_{21} e_1 = s_{21}, \; s_{12} e_1 = e_1 s_{21} = 0 \]
Any $S$-module of dimension $n$ is therefore given by matrices $E_1,
S_{12}, S_{21} \in M_n(k)$ satisfying the matric equations
    \[ S_{12}^2 = S_{21}^2 = 0, \; E_1^2 = E_1, \; E_1 S_{12} = S_{12},
    \; S_{21} E_1 = S_{21}, \; S_{12} E_1 = E_1 S_{21} = 0 \]
The $S$-modules represented by $(E_1, S_{12}, S_{21})$ and $(E'_1,
S'_{12}, S'_{21})$ are isomorphic if and only if there is an invertible
matrix $G \in M_n(k)$ such that $G E_1 G^{-1} = E_1', \; G S_{12}
G^{-1} = S_{12}', G S_{21} G^{-1} = S_{21}'$. Using this
characterization, it is a straight-forward but tedious task to classify
all $S$-modules of dimension $n$ up to isomorphism for a given integer
$n \ge 1$.

Let us first remark that for any $S$-module of dimension $n = 1$,
$\rho$ factorizes through the commutativization $k^2$ of $S$. It
follows that there are exactly two non-isomorphic simple $S$-modules of
dimension one, $T_{1,1}$ and $T_{1,2}$, and the corresponding
deformations of $\sfam$ are
    \[ M_{1,i} = \sfam_S \otimes_S T_{1,i} \cong M_i \quad \text{for }
    i = 1,2 \]
This reflects that $M_1$ and $M_2$ are rigid as left $D$-modules.

We obtain the following list of $S$-modules of dimension $n = 2$, up to
isomorphism. We have used that, without loss of generality, we may
assume that $E_1$ has Jordan form:
\begin{align}
E_1 &= \left( \begin{matrix} 0 & 0 \\ 0 & 0 \end{matrix} \right) &
S_{12} &= \left( \begin{matrix} 0 & 0 \\ 0 & 0 \end{matrix} \right) &
S_{21} &= \left( \begin{matrix} 0 & 0 \\ 0 & 0 \end{matrix} \right) \\
E_1 &= \left( \begin{matrix} 1 & 0 \\ 0 & 1 \end{matrix} \right) &
S_{12} &= \left( \begin{matrix} 0 & 0 \\ 0 & 0 \end{matrix} \right) &
S_{21} &= \left( \begin{matrix} 0 & 0 \\ 0 & 0 \end{matrix} \right) \\
E_1 &= \left( \begin{matrix} 1 & 0 \\ 0 & 0 \end{matrix} \right) &
S_{12} &= \left( \begin{matrix} 0 & 0 \\ 0 & 0 \end{matrix} \right) &
S_{21} &= \left( \begin{matrix} 0 & 0 \\ 0 & 0 \end{matrix} \right) \\
E_1 &= \left( \begin{matrix} 1 & 0 \\ 0 & 0 \end{matrix} \right) &
S_{12} &= \left( \begin{matrix} 0 & 0 \\ 0 & 0 \end{matrix} \right) &
S_{21} &= \left( \begin{matrix} 0 & 0 \\ 1 & 0 \end{matrix} \right) \\
E_1 &= \left( \begin{matrix} 1 & 0 \\ 0 & 0 \end{matrix} \right) &
S_{12} &= \left( \begin{matrix} 0 & 1 \\ 0 & 0 \end{matrix} \right) &
S_{21} &= \left( \begin{matrix} 0 & 0 \\ 0 & 0 \end{matrix} \right) \\
E_1 &= \left( \begin{matrix} 1 & 0 \\ 0 & 0 \end{matrix} \right) &
S_{12} &= \left( \begin{matrix} 0 & 1 \\ 0 & 0 \end{matrix} \right) &
S_{21} &= \left( \begin{matrix} 0 & 0 \\ a & 0 \end{matrix} \right)
\quad \text{ for } a \in k^*
\end{align}
We shall write $T_{2,1}$ -- $T_{2,5}$ and $T_{2,6,a}$ for the
corresponding $S$-modules of dimension two. Notice that $T_{2,6,a}$ is
simple for all $a \in k^*$, while $T_{2,1}$ -- $T_{2,5}$ are extensions
of simple $S$-modules of dimension one. The corresponding deformations
of $\sfam$ are given by
    \[ M_{2,6,a} = \sfam_S \otimes_S T_{2,6,a} \quad \text{for } a \in
    k^* \]
In fact, one may show that $M_{2,6,a} \cong D/D \cdot (t \partial - a)$
for any $a \in k^*$. In particular, $M_{2,6,a}$ is a simple $D$-module
if $a \not \in \z$, and in this case $M_{2,6,a} \cong M_{2,6,b}$ if and
only if $a-b \in \z$. Furthermore, $M_{2,6,-1} \cong D/D \cdot \partial
\; t$, $M_{2,6,n} \cong M_1$ for $n = 1,2,\dots$, and $M_{2,6,-n} \cong
M_2$ for $n = 2,3,\dots$.

We obtain the following list of $S$-modules of dimension $n = 3$, up to
isomorphism. We have used that, without loss of generality, we may
assume that $E_1$ has Jordan form:
\begin{align} \setcounter{equation}{0}
E_1 &= \left( \begin{matrix} 0 & 0 & 0 \\ 0 & 0 & 0 \\ 0 & 0 & 0
\end{matrix} \right) &
S_{12} &= \left( \begin{matrix} 0 & 0 & 0 \\ 0 & 0 & 0 \\ 0 & 0 & 0
\end{matrix} \right) &
S_{21} &= \left( \begin{matrix} 0 & 0 & 0 \\ 0 & 0 & 0 \\ 0 & 0 & 0
\end{matrix} \right) \\
E_1 &= \left( \begin{matrix} 1 & 0 & 0 \\ 0 & 1 & 0 \\ 0 & 0 & 1
\end{matrix} \right) &
S_{12} &= \left( \begin{matrix} 0 & 0 & 0 \\ 0 & 0 & 0 \\ 0 & 0 & 0
\end{matrix} \right) &
S_{21} &= \left( \begin{matrix} 0 & 0 & 0 \\ 0 & 0 & 0 \\ 0 & 0 & 0
\end{matrix} \right) \\
E_1 &= \left( \begin{matrix} 1 & 0 & 0 \\ 0 & 0 & 0 \\ 0 & 0 & 0
\end{matrix} \right) &
S_{12} &= \left( \begin{matrix} 0 & 0 & 0 \\ 0 & 0 & 0 \\ 0 & 0 & 0
\end{matrix} \right) &
S_{21} &= \left( \begin{matrix} 0 & 0 & 0 \\ 0 & 0 & 0 \\ 0 & 0 & 0
\end{matrix} \right) \\
E_1 &= \left( \begin{matrix} 1 & 0 & 0 \\ 0 & 0 & 0 \\ 0 & 0 & 0
\end{matrix} \right) &
S_{12} &= \left( \begin{matrix} 0 & 0 & 0 \\ 0 & 0 & 0 \\ 0 & 0 & 0
\end{matrix} \right) &
S_{21} &= \left( \begin{matrix} 0 & 0 & 0 \\ 1 & 0 & 0 \\ 0 & 0 & 0
\end{matrix} \right) \\
E_1 &= \left( \begin{matrix} 1 & 0 & 0 \\ 0 & 0 & 0 \\ 0 & 0 & 0
\end{matrix} \right) &
S_{12} &= \left( \begin{matrix} 0 & 1 & 0 \\ 0 & 0 & 0 \\ 0 & 0 & 0
\end{matrix} \right) &
S_{21} &= \left( \begin{matrix} 0 & 0 & 0 \\ 0 & 0 & 0 \\ 0 & 0 & 0
\end{matrix} \right) \\
E_1 &= \left( \begin{matrix} 1 & 0 & 0 \\ 0 & 0 & 0 \\ 0 & 0 & 0
\end{matrix} \right) &
S_{12} &= \left( \begin{matrix} 0 & 1 & 0 \\ 0 & 0 & 0 \\ 0 & 0 & 0
\end{matrix} \right) &
S_{21} &= \left( \begin{matrix} 0 & 0 & 0 \\ 0 & 0 & 0 \\ 1 & 0 & 0
\end{matrix} \right) \\
E_1 &= \left( \begin{matrix} 1 & 0 & 0 \\ 0 & 0 & 0 \\ 0 & 0 & 0
\end{matrix} \right) &
S_{12} &= \left( \begin{matrix} 0 & 1 & 0 \\ 0 & 0 & 0 \\ 0 & 0 & 0
\end{matrix} \right) &
S_{21} &= \left( \begin{matrix} 0 & 0 & 0 \\ b & 0 & 0 \\ 0 & 0 & 0
\end{matrix} \right)
\quad \text{ for } b \in k^* \\
E_1 &= \left( \begin{matrix} 1 & 0 & 0 \\ 0 & 1 & 0 \\ 0 & 0 & 0
\end{matrix} \right) &
S_{12} &= \left( \begin{matrix} 0 & 0 & 0 \\ 0 & 0 & 0 \\ 0 & 0 & 0
\end{matrix} \right) &
S_{21} &= \left( \begin{matrix} 0 & 0 & 0 \\ 0 & 0 & 0 \\ 0 & 0 & 0
\end{matrix} \right) \\
E_1 &= \left( \begin{matrix} 1 & 0 & 0 \\ 0 & 1 & 0 \\ 0 & 0 & 0
\end{matrix} \right) &
S_{12} &= \left( \begin{matrix} 0 & 0 & 0 \\ 0 & 0 & 0 \\ 0 & 0 & 0
\end{matrix} \right) &
S_{21} &= \left( \begin{matrix} 0 & 0 & 0 \\ 0 & 0 & 0 \\ 0 & 1 & 0
\end{matrix} \right) \\
E_1 &= \left( \begin{matrix} 1 & 0 & 0 \\ 0 & 1 & 0 \\ 0 & 0 & 0
\end{matrix} \right) &
S_{12} &= \left( \begin{matrix} 0 & 0 & 0 \\ 0 & 0 & 1 \\ 0 & 0 & 0
\end{matrix} \right) &
S_{21} &= \left( \begin{matrix} 0 & 0 & 0 \\ 0 & 0 & 0 \\ 0 & 0 & 0
\end{matrix} \right) \\
E_1 &= \left( \begin{matrix} 1 & 0 & 0 \\ 0 & 1 & 0 \\ 0 & 0 & 0
\end{matrix} \right) &
S_{12} &= \left( \begin{matrix} 0 & 0 & 0 \\ 0 & 0 & 1 \\ 0 & 0 & 0
\end{matrix} \right) &
S_{21} &= \left( \begin{matrix} 0 & 0 & 0 \\ 0 & 0 & 0 \\ 1 & 0 & 0
\end{matrix} \right) \\
E_1 &= \left( \begin{matrix} 1 & 0 & 0 \\ 0 & 1 & 0 \\ 0 & 0 & 0
\end{matrix} \right) &
S_{12} &= \left( \begin{matrix} 0 & 0 & 0 \\ 0 & 0 & 1 \\ 0 & 0 & 0
\end{matrix} \right) &
S_{21} &= \left( \begin{matrix} 0 & 0 & 0 \\ 0 & 0 & 0 \\ 0 & c & 0
\end{matrix} \right)
\quad \text{ for } c \in k^*
\end{align}
We shall write $T_{3,1}$ -- $T_{3,6}$, $T_{3,7,b}$, $T_{3,8}$ --
$T_{3,11}$, and $T_{3,12,c}$ for the corresponding $S$-modules of
dimension three. Notice that all $S$-modules of dimension three are
extensions of simple $S$-modules of dimension one and two, so there are
no simple $S$-modules of dimension $n = 3$.

We remark that if $T$ is a simple $S$-module, then $\rho: S \to
\enm_k(T)$ is a surjective ring homomorphism. Hence it seems unlikely
that there are any simple $S$-modules of dimension $n > 3$.

Finally, we remark that the commutative deformation functor $\defm{M}:
\lcat \to \sets$ of the direct sum $M = M_1 \oplus M_2$ has
pro-representing hull $(H = k[[ s_{12}, s_{21} ]], M_H)$, and an
algebraization $(S = k[ s_{12}, s_{21} ],M_S)$. It is not difficult to
find the family $M_S$ in this case. In fact, for any point $(\alpha,
\beta) \in \spec S = \aff^2_k$, the left $D$-module $M_{\alpha, \beta}
= M_S \otimes_S S/( s_{12} - \alpha, s_{21} - \beta )$ is given by
\begin{align*}
    M_{0,0} &\cong M_1 \oplus M_2 & \\
    M_{\alpha,0} &\cong D/D \cdot (\partial \; t) & \text{for } \alpha
    \neq 0 \\
    M_{\alpha,\beta} &\cong D/D \cdot ( t \; \partial - \alpha \beta )
    & \text{for } \beta \neq 0
\end{align*}
We see that we obtain exactly the same isomorphism classes of left
$D$-modules as commutative deformations of $M = M_1 \oplus M_2$ as we
obtained as noncommutative deformations of the family $\sfam = \{ M_1,
M_2 \}$. However, the points of the pro-representing hull
    \[ H = \left( \begin{matrix} H_{11} & H_{12} \\ H_{21} & H_{22}
    \end{matrix} \right) \cong \left( \begin{matrix} k[[ s_{12} s_{21}
    ]] & \langle s_{12} \rangle \\ \langle s_{21} \rangle & k[[ s_{21}
    s_{12} ]] \end{matrix} \right) \]
of the noncommutative deformation functor $\defs$ give a much better
geometric picture of the local structure of the moduli space of left
$D$-modules.

\bibliographystyle{plain}
\bibliography{defm-gb}

\end{document}